\documentclass[a4paper,12pt]{amsart}
\usepackage{amssymb,amsthm}
\usepackage{amsmath}
\usepackage{tikz}
\usepackage{graphicx}
\usepackage{subcaption}
\usepackage[shortlabels]{enumitem}
\usetikzlibrary{arrows}
\usepackage{float}
\usepackage{graphicx}
\usepackage{subcaption}
\graphicspath{{ToricFigs/}}
\usetikzlibrary{decorations.pathmorphing}
\tikzset{snake it/.style={decorate, decoration=snake}}

\makeatletter
\newcommand*{\rom}[1]{\expandafter\@slowromancap\romannumeral #1@}
\makeatother

\numberwithin{equation}{section}

\theoremstyle{plain}
\newtheorem{theorem}{Theorem}
\numberwithin{theorem}{section}
\newtheorem{proposition}[theorem]{Proposition}

\theoremstyle{definition}
\newtheorem{definition}[theorem]{Definition}
\theoremstyle{remark}
\newtheorem{remark}[theorem]{Remark}
\theoremstyle{remark}
\newtheorem{discussion}[theorem]{Discussion}

\DeclareMathOperator{\sinc}{sinc}

\usepackage{fullpage}

\newcommand{\be}{\begin{equation}}
\newcommand{\ee}{\end{equation}}
\newcommand{\bea}{\begin{eqnarray}}
\newcommand{\eea}{\end{eqnarray}}
\newcommand{\bean}{\begin{eqnarray*}}
\newcommand{\eean}{\end{eqnarray*}}

\newcommand{\bel}[1]{\begin{equation}\label{#1}}

\newcommand{\eel}[1]{{\label{#1}\end{equation}}}

\usepackage{datetime}

\title[Compton scattering tomography in translational geometries]{Compton scattering tomography in translational geometries\\{\footnotesize\ddmmyyyydate\today~\currenttime}}

\author{James Webber}
\address{Department of Electrical and Computer
Engineering, Tufts University, Medford, MA USA}
\email{James.Webber@tufts.edu}
\author{Eric L. Miller}
\address{Department of Electrical and Computer
Engineering, Tufts University, Medford, MA USA}
\email{elmiller@ece.tufts.edu}

\begin{document}

\begin{abstract}
Here we present new $L^2$ injectivity results for 2-D and 3-D Compton scattering tomography (CST) problems in translational geometries. The results are proven through the explicit inversion of a new toric section and apple Radon transform, which describe novel 2-D and 3-D acquisition geometries in CST. The geometry considered has potential applications in airport baggage screening and threat detection. We also present a generalization of our injectivity results in 3-D to Radon transforms which describe the integrals of the charge density over the surfaces of revolution of a class of $C^1$ curves.
\end{abstract}

\maketitle

\section{Introduction} 
\begin{figure}[!h]
\centering
\begin{tikzpicture}[scale=6]
\draw [very thick] (-1,0)--(1,0)node[right] {$\{z=2-r_m\}$};
\draw [very thick] (-1,1)--(1,1)node[right] {$d$};
\draw [<->] (-0.1,1)--(-0.1,1.25);
\node at (-0.14,1.1125) {$1$};
\draw [<->] (-0.1,1.25)--(-0.1,1.5);
\node at (-0.14,1.3725) {$1$};
\draw [<->] (-0.1,0.75)--(-0.1,1);
\node at (-0.14,0.8725) {$1$};
\draw [very thick] (-1,1.5)--(1,1.5)node[right] {$s$};
\draw [->] (0,0.75)--(0,1.7);
\draw [->] (0,0.75)--(1,0.75);
\draw (0.6,1.25) circle (0.65);
\draw (-0.6,1.25) circle (0.65);
\node at (-0.6,2) {$T(r)$ (2-D), $A(r)$ (3-D)};
\draw [<->] (0.5+0.1,1.25)--(0,1);
\node at (0.24+0.08,1.1) {$r$};
\draw [<->] (0.5+0.1,1.25)--(0,1.25);
\node at (0.24+0.05,1.3) {$R$};
\draw [<->] (0.5+0.1,0.6)--(0.5+0.1,0.75);
\node at (0.55,0.67) {$|h|$};
\node at (0.04,1.75) {$z$};
\node at (1.04,0.45+0.25) {$x$};
\fill[blue!40!white] (-0.5,0.2) rectangle (0,0.4);
\fill[purple!40!white] (-0.75,0.5) circle (0.25);
\draw [->] (-1.1,0.1)--(-0.75,0.25);
\node at (-1.1,0.14) {$f$};
\node at (-0.03,0.47+0.25) {$O$};
\end{tikzpicture}
\caption{Vertical torus/toric section scanning geometry. We denote $T(r)$ to be the union of the two circles drawn above, and the apple $A(r)$ to be the surface of revolution of $T(r)$ about $z$. $s$ and $d$ label the source and detector lines respectively, $r$ is the radius of the torus, $R=\sqrt{r^2-1}$ is the distance from the center of the torus to the centre of the torus tube and $h=2-r$.}
\label{fig1}
\end{figure}
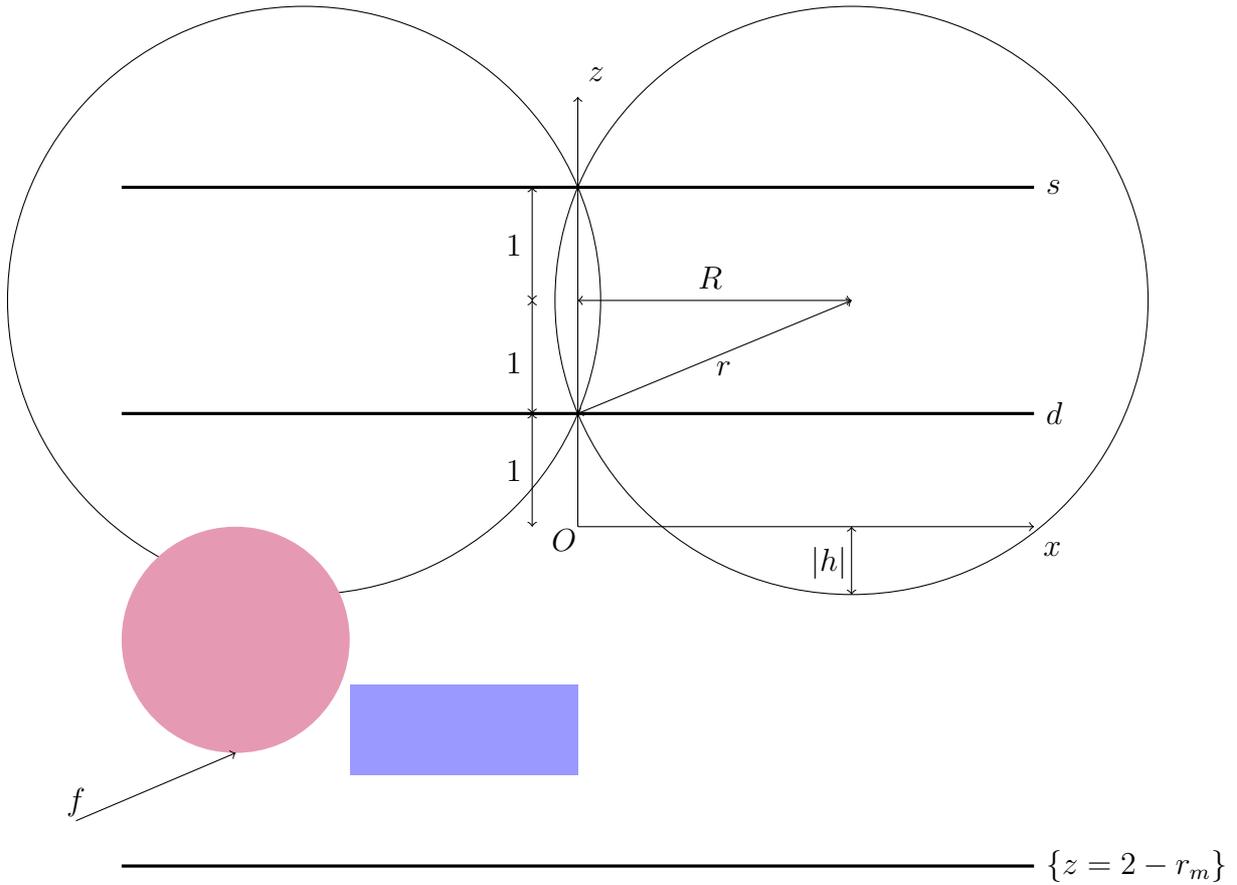

Here we present novel inversion formulae for two and three dimensional scanning modalities in CST, which have potential applications in airport baggage screening and threat detection. We consider the scanning geometry illustrated in figure \ref{fig1}. The electron charge density $f$ (represented by a real valued function) is translated in the $y$ direction (on a conveyor belt) and illuminated by a line of monochromatic photon sources. The scattered intensity is then collected by the detector array. The source and detector array are assumed to be slightly offset (similarly to the machine geometry considered in \cite{will}) in the $y$ direction so that the detectors do not block any incoming photons. Given the small offset we model the sources and detectors to lie in the same ($xz$) plane, as in \cite{me}, for elegance of mathematical derivation.

The works of \cite{2D1,2D2,2D3,2D4,2D5,2D6,2D7,2D8} (2-D geometries) and \cite{3D1,3D2,3D3,3D4,3D5} (3-D geometries) respectively consider CST problems in a variety of two and three dimensional scanning geometries. In electron density reconstruction in CST, the scattered intensity is modelled as integrals of the electron density over circular arcs (in 2-D with collimated detectors \cite{2D2}), toric sections (in 2-D with uncollimated detectors \cite{me,2D8}) and spindle tori (in 3-D with no collimation \cite{3D2,3D3}). In gamma ray source reconstruction using CST, the scattered intensity is modelled as integrals of the source intensity over V-lines (or broken rays, in 2-D \cite{2D4,2D6}) and cones (in 3-D \cite{3D1,3D4}). In contrast to these efforts, here we consider the problem of reconstructing the electron density from a set of vertical (with axis of revolution parallel to the $z$ axis) toric sections (2-D) and apples (3-D, the apple is the part of the spindle torus which corresponds to backscattered photons) which are translated along the line (2-D) and in the plane (3-D) respectively. The nature of the data acquisition is such that the previous results from the literature are insufficient so as to provide an explicit solution. We aim to provide such explicit inversion results here in this paper. 

The inversion process we present relies on Paley-Weiner-Schwartz ideas (using analytic continuation as in \cite{will}) and the explicit inversion of 1-D Volterra operators in the Fourier domain. In \cite{3D1}, the source intensity is reconstructed from a set of vertical cones translated on the plane (the data dimensions are the cone opening angle, and a 2-D translation). Here, after the data is transformed to the Fourier space, the inversion is carried out by inverse Hankel transform of a set of 1-D integral equations with Bessel kernels. In section \ref{3D} we also discover Bessel function kernels in the Fourier domain, but in our case we use the theory of Volterra integral operators \cite{Tric} to invert the resulting set of 1-D integral equations. In 3-D our geometry is akin to configuration (d) of \cite[page 5]{3D2}, but we only consider the vertical spindle tori. In \cite{3D2} they provide microlocal inversion results (modulo smoothing) for the proposed geometry. We provide injectivity results and exact explicit inversion formulae for geometry (d) in this paper. In \cite{2D8} injectivity results are proven for a toric section transform with rotational invariance. Here, after an expansion to the Fourier series, the solution was obtained using the theory of Cormack \cite{cor}. We consider a toric section transform with translational invariance, whereby our solution is obtained using the theory of \cite{Tric}. 

In section \ref{pre} we state some definitions and preliminary results that will be needed to prove our main theorems. In section \ref{2D} we explain the parameterization of the toric section curves in the translation geometry and introduce a new toric section transform $\mathcal{T}$. We then prove the boundedness of $\mathcal{T}$ in $L^2$ (this is a typical property of Radon transforms as smoothing operators) before going on to prove our first main theorem which proves the injectivity and explicit invertibility of $\mathcal{T}$ on the domain of compactly supported $L^2$ functions.

In section \ref{3D} we consider the 3-D case, where we introduce a new apple Radon transform $\mathcal{A}$. As in the 2-D case in section \ref{2D} we prove the continuity of $\mathcal{A}$ in $L^2$ before going on to prove our second main theorem, which shows the injectivity and explicit invertibility of $\mathcal{A}$ on the $L^2$ domain. The proof uses similar ideas to the 2-D case, except in 3-D the resulting kernels are the finite sum of Bessel functions, whereas in 2-D we see Volterra operators with cosine kernels.

In section \ref{gen} we introduce the generalized Radon transform $R_A$, which describes the integrals of a density over the surfaces of revolution of $C^1$ curves, in the translational geometry. Example $C^1$ curves would include the previously considered semicircles (for apples) and straight lines (for cones \cite{3D1}). We show the injectivity of $R_A$ on the domain of compactly supported continuous functions in this case. Our proof uses the ideas of Cormack \cite{cor}, which require continuity in the target function. 
\section{Preliminary results and definitions}
\label{pre}
Here we state some definitions and preliminary results which will be used in our theorems. 
\begin{definition}
Let $f\in L^2_0(\mathbb{R}^n)$. Then we define the Fourier transform $\hat{f}$ of $f$ in terms of angular frequency
\begin{equation}
\hat{f}(\omega)=(2\pi)^{-n/2}\int_{\mathbb{R}^n}f(x)e^{-ix\cdot \omega}\mathrm{d}x.
\end{equation}
\end{definition}
We state the Plancharel theorem \cite{plan}.
\begin{theorem}[Plancharel theorem]
Let $f\in L^2_0(\mathbb{R}^n)$. Then $\hat{f}\in L^2(\mathbb{R}^n)$ and
\begin{equation}
\|f\|_{L^2_0(\mathbb{R}^n)}=\|\hat{f}\|_{L^2(\mathbb{R}^n)}.
\end{equation}
\end{theorem}
From \cite[page 22]{hor1}, we have the Paley-Weiner-Schwartz theorem.
\begin{theorem}[Paley-Weiner-Schwartz]
Let $\mathcal{E}'(\mathbb{R}^n)$ be the set of distributions of compact support in $\mathbb{R}^n$ and let $f\in\mathcal{E}'(\mathbb{R}^n)$. Let $U=\text{supp}(f)$ and let $H_U$ be the support function of $U$, defined as
\begin{equation}
H_U(\omega)=\sup_{x\in U}(x\cdot \omega).
\end{equation}
Then the Fourier transform $\hat{f}$ is an entire analytic function. Additionally, there exists a constant $c>0$ and integer $N$ such that
\begin{equation}
|\hat{f}(\omega)|\leq c(1+|\omega|)^Ne^{H_U(\text{Im}(\omega))},\ \ \ \forall \omega\in\mathbb{C}^n,
\end{equation}
where $\text{Im}(\omega)$ denotes the imaginary part of $\omega$. Conversely, let $F$ be an entire analytic function, let $H$ be the support function of a compact convex set $K\subset \mathbb{R}^n$ and suppose that
\begin{equation}
|F(\omega)|\leq c(1+|\omega|)^Ne^{H(\text{Im}(\omega))},\ \ \ \forall\omega\in\mathbb{C}^n
\end{equation}
for some constants $c$ and $N$. Then there exists a unique $f\in\mathcal{E}'(\mathbb{R}^n)$ such that $\hat{f}=F$ and $\text{supp}(f)\subset K$.
\end{theorem}
We now state some results on Volterra type integral equations from \cite[page 10]{Tric}.
\begin{definition}
We define a Volterra equation of the second kind to be an equation of the form
\begin{equation}
g(x)=\lambda\int_0^xK(x,y)f(y)\mathrm{d}y+f(x)
\end{equation}
with real valued kernel $K$ on a triangle $T'=\{0<x<x', 0<y<x\}$. $K$ is said to be an $L^2$ kernel if
\begin{equation}
\|K\|^2_{L^2(T')}=\int_0^{x'}\int_0^x K^2(x,y)\mathrm{d}y\mathrm{d}x\leq N^2
\end{equation}
for some $N>0$.
\end{definition}
\begin{theorem}
\label{volt}
Let $f\in L^2([0,x'])$ and let $K(x,y)$ be an $L^2$ kernel on $T'=\{0<x<x', 0<y<x\}$ for some $x'>0$. Then the Volterra integral equation of the second kind
\begin{equation}
g(x)=\lambda\int_0^xK(x,y)f(y)\mathrm{d}y+f(x)
\end{equation}
has one and only one solution in $L^2([0,x'])$, and the solution is given by the formula
\begin{equation}
f(x)=\lambda\int_0^xH(x,y;\lambda)g(y)\mathrm{d}y+g(y),
\end{equation}
where
\begin{equation}
H(x,y;\lambda)=\sum_{\nu=0}^{\infty}\lambda^{\nu}K_{\nu+1}(x,y)
\end{equation}
and the iterated kernels $K_{\nu}$ are defined by
$$K_1(x,y)=K(x,y)$$
and
$$K_{\nu+1}(x,y)=\int_0^xK(x,z)K_{\nu}(z,y)\mathrm{d}z$$
for $\nu\geq 1$.
\end{theorem}
Next we have Young's inequality \cite[Theorem 0.3.1]{young}
\begin{theorem}[Young's inequality]
Let $X$ and $Y$ be measurable spaces and let $q,p,t\in\mathbb{Z}^+\backslash\{0\}$ be such that $\frac{1}{q}=\frac{1}{p}+\frac{1}{t}-1$. If $K: X\times Y\to \mathbb{R}$ is such that there exists a $C>0$ with
\begin{equation}
\int_{X}|K(x,y)|^t\mathrm{d}x\leq C^t
\end{equation}
and
\begin{equation}
\int_{Y}|K(x,y)|^t\mathrm{d}y\leq C^t.
\end{equation}
Then for an integrable function $f$
\begin{equation}
\int_{X}\left|\int_YK(x,y)f(y)\mathrm{d}y\right|^q\mathrm{d}x\leq C^q\left(\int_Y|f(y)|^p\mathrm{d}y\right)^{\frac{q}{p}}.
\end{equation}
\end{theorem}
\section{Two dimensional Compton tomography under translations}
\label{2D}
Let $f\in L^2_0(\Omega)$ be the electron charge density, where $\Omega\subset\{2-r_m<z<1\}$ for some $r_m>1$, is compact. Here $r_m$ controls the (fixed) depth of the scanning tunnel as is depicted in figure \ref{fig1}. Let $r=2-h$ and $R=\sqrt{r^2-1}$ (as is illustrated in figure \ref{fig1}). Let $S_1(r),\ldots,S_4(r)$ be the four semicircles, radius $r$, whose disjoint union $\cup_{j=1}^4S_j(r)=T(r)$ is the toric section in figure \ref{fig1}, and such that the parameterization of each $S_j(r)$ for $j=1,2,3,4$ is given by the vertical axis coordinate
\begin{equation}
\begin{split}
x_1&=R+\sqrt{r^2-(z-2)^2}, \ \ \ x_2=R-\sqrt{r^2-(z-2)^2}\ \ \ \text{(right hand circle)}\\
x_3&=-R+\sqrt{r^2-(z-2)^2}, \ \ \ x_4=-R-\sqrt{r^2-(z-2)^2},\ \ \ \text{(left hand circle)}
\end{split}
\end{equation}
for $2-r<z<1$. We define the toric section transform $\mathcal{T}: L^2_0(\Omega)\to \mathcal{T}(L^2_0(\Omega))$ in the translational geometry as
\begin{equation}
\begin{split}
\mathcal{T}f(x_0,r)=\int_{T_{x_0}(T(r))}f\mathrm{d}s=\sum_{j=1}^4\int_{T_{x_0}(S_j(r))}f\mathrm{d}s_j,
\end{split}
\end{equation}
where $\mathrm{d}s_j$ denotes the arc length measure on $S_j$ and $T_{x_0}(x,z)=(x+x_0,z)$ denotes a translation along the $x$ axis of length $x_0$. 
\begin{proposition}
Let $f\in L^2_0(\Omega)$, and let $f_1$ be defined as $f_1(x,z)=f(x,2-z)$. Then
\begin{equation}
\begin{split}
\mathcal{T}f(x_0,r)=\int_1^r\frac{r}{\sqrt{r^2-z^2}}&\sum_{j=1}^2f_1(R+(-1)^{j}\sqrt{r^2-z^2}+x_0,z)\\
&+f_1(-R+(-1)^{j}\sqrt{r^2-z^2}+x_0,z)\mathrm{d}z.
\end{split}
\end{equation}
\begin{proof}
The arc length measure is given by
\begin{equation}
\label{arc}
\mathrm{d}s_j=\mathrm{d}z\sqrt{1+\left(\frac{\mathrm{d}x_j}{\mathrm{d}z}\right)^2}=\frac{r}{\sqrt{r^2-(z-2)^2}}\mathrm{d}z.
\end{equation}
From which it follows that
\begin{equation}
\label{Tform}
\begin{split}
\mathcal{T}f(x_0,r)&=\sum_{j=1}^4\int_{T_{x_0}(S_j(r))}f\mathrm{d}s_j\\
&=\int_h^1\frac{r}{\sqrt{r^2-(z-2)^2}}\sum_{j=1}^4f(x_j+x_0,z)\mathrm{d}z \\
&=\int_1^r\frac{r}{\sqrt{r^2-z^2}}\sum_{j=1}^2f_1(R+(-1)^{j}\sqrt{r^2-z^2}+x_0,z)\\
&+f_1(-R+(-1)^{j}\sqrt{r^2-z^2}+x_0,z)\mathrm{d}z,
\end{split}
\end{equation}
which completes the proof.
\end{proof}
\end{proposition}
We have the result which proves the boundedness of $\mathcal{T}$ in $L^2_0$.
\begin{proposition}
\label{prp2}
Let $r_m>1$ and let $\Omega\subset\{2-r_m<z<1\}$ be compact. Then the image $\mathcal{T}(L^2_0(\Omega))\subset L^2([1,r_m]\times\mathbb{R})$.
\begin{proof}
Taking the Fourier transform in the $x_0$ variable of \eqref{Tform} yields
\begin{equation}
\begin{split}
\widehat{\mathcal{T}f}(\omega_1,r)&=\int_{\mathbb{R}}\mathcal{T}f(x_0,r)e^{-ix_0\omega_1}\mathrm{d}x_0\\
&=\int_1^r\frac{rK(r,z)}{\sqrt{r^2-z^2}}\hat{f_1}(\omega_1,z)\mathrm{d}z,
\end{split}
\end{equation}
where $\omega_1$ is dual to $x$ (or $x_0$) and
\begin{equation}
\hspace{-1cm}
\begin{split}
K(r,z)&=\sum_{j=1}^2\exp\left( -i\omega_1(R+(-1)^{j}\sqrt{r^2-z^2})\right)+\exp\left(-i\omega_1(-R+(-1)^{j}\sqrt{r^2-z^2})\right),
\end{split}
\end{equation}
where $|K|\leq 4$ is bounded. We have
\begin{equation}
\begin{split}
\|\mathcal{T}f\|^2_{L^2([1,r_m]\times\mathbb{R})}&=\|\widehat{\mathcal{T}f}\|^2_{L^2([1,r_m]\times\mathbb{R})}\ \ \ \text{by the Plancharel theorem}\\
&=\int_{-\infty}^{\infty}\int_1^{r_m}\left(\int_1^r\frac{rK(r,z)}{\sqrt{r^2-z^2}}\hat{f_1}(\omega_1,z)\mathrm{d}z\right)^2\mathrm{d}r\mathrm{d}\omega_1\\
&\leq 16r^2_m\int_{-\infty}^{\infty}\int_1^{r_m}\left(\int_1^r\frac{\hat{f_1}(\omega_1,z)}{\sqrt{r^2-z^2}}\mathrm{d}z\right)^2\mathrm{d}r\mathrm{d}\omega_1\\
&\leq\frac{(4r_m)^2}{2}\int_{-\infty}^{\infty}\int_1^{r_m}\left(\int_1^r\hat{f_1}(\omega_1,z)L(r,z)\mathrm{d}z\right)^2\mathrm{d}r\mathrm{d}\omega_1\\,
\end{split}
\end{equation}
where $L(r,z)=(r-z)^{-1/2}$ is defined on the triangle $H=\{1<r<r_m, 1<z<r\}$. Let us extend the domain of $L$ to $[1,r_m]^2$, in the sense that $L(r,z)=0$ on $[1,r_m]^2\backslash H$. Let $t=1$. Then we have
\begin{equation}
\int_{[1,r_m]}|L(r,z)|^t\mathrm{d}r=\int_z^{r_m}(r-z)^{-1/2}\mathrm{d}r=2\sqrt{r_m-z}\leq 2\sqrt{r_m-1}=C^t
\end{equation}
and
\begin{equation}
\int_{[1,r_m]}|L(r,z)|^t\mathrm{d}z=\int_1^r(r-z)^{-1/2}\mathrm{d}z=2\sqrt{r-1}\leq 2\sqrt{r_m-1}=C^t.
\end{equation}
Setting $p=q=2$ and $X=Y=[1,r_m]$, it follows from Young's inequality that
\begin{equation}
\begin{split}
\int_1^{r_m}\left(\int_1^r\hat{f_1}(\omega_1,z)L(r,z)\mathrm{d}z\right)^2\mathrm{d}r&=\int_{[1,r_m]}\left(\int_{[1,r_m]}\hat{f_1}(\omega_1,z)L(r,z)\mathrm{d}z\right)^2\mathrm{d}r\\
&\leq C^2\int_{[1,r_m]}|\hat{f_1}(\omega_1,z)|^2\mathrm{d}z.
\end{split}
\end{equation}
Hence we have
\begin{equation}
\begin{split}
\|\mathcal{T}f\|^2_{L^2([1,r_m]\times\mathbb{R})}&\leq\frac{(4r_m)^2}{2}\int_{-\infty}^{\infty}\int_1^{r_m}\left(\int_1^r\hat{f_1}(\omega_1,z)L(r,z)\mathrm{d}z\right)^2\mathrm{d}r\mathrm{d}\omega_1\\
&\leq 8(r_mC)^2\int_{-\infty}^{\infty}\|\hat{f_1}(\cdot,\omega_1)\|^2_{L^2([1,r_m])}\mathrm{d}\omega_1\\
&=8(r_mC)^2\|f\|^2_{L^2(\mathbb{R}^2)}<\infty,
\end{split}
\end{equation}
which completes the proof.
\end{proof}
\end{proposition}
We now have our first main theorem
\begin{theorem}
\label{main1}
Let $r_m>1$ and let $\Omega\subset\{2-r_m<z<1\}$ be compact. Then the toric section transform $\mathcal{T}: L^2_0(\Omega)\to L^2([1,r_m]\times\mathbb{R})$ is injective.
\begin{proof}
From Proposition \ref{prp2} we have the representation of the toric section transform in the Fourier domain
\begin{equation}
\begin{split}
\widehat{\mathcal{T}f}(\omega_1,r)=\int_1^r\frac{rK(r,z)}{\sqrt{r^2-z^2}}\hat{f_1}(\omega_1,z)\mathrm{d}z,
\end{split}
\end{equation}
where we now have the simplified form for $K$
\begin{equation}
\hspace{-1cm}
\begin{split}
K(r,z)&=\sum_{j=1}^2\exp\left( -i\omega_1(R+(-1)^{j}\sqrt{r^2-z^2})\right)+\exp\left(-i\omega_1(-R+(-1)^{j}\sqrt{r^2-z^2})\right)\\
&=\left(\exp( -i\omega_1 R)+\exp(i\omega_1 R)\right)\sum_{j=1}^2\exp\left( -i\omega_1(-1)^{j}\sqrt{r^2-z^2}\right)\\
&=4\cos(\omega_1 R)\cos(\omega_1\sqrt{r^2-z^2}).
\end{split}
\end{equation}
After making the substitution $z=z^2$ and letting $\hat{f_2}(\omega_1,z)=\frac{\hat{f_1}(\omega_1,\sqrt{z})}{2\sqrt{z}}$, we have
\begin{equation}
\begin{split}
\widehat{\mathcal{T}_1f}(\omega_1,r)&=\frac{\widehat{\mathcal{T}f}(\omega_1,\sqrt{r})}{4\sqrt{r}}\\
&=\cos(\omega_1 \sqrt{r-1})\int_1^r\frac{\cos(\omega_1\sqrt{r-z})}{\sqrt{r-z}}\hat{f_2}(\omega_1,z)\mathrm{d}z.
\end{split}
\end{equation}
For $|\omega_1|<\frac{\pi}{2\sqrt{r^2_m-1}}$ and $r<r^2_m$, $\cos(\omega_1\sqrt{r-1})>0$ and we have
\begin{equation}
\begin{split}
\hat{g}(\omega_1,r)&=\frac{\widehat{\mathcal{T}_1f}(\omega_1,r)}{\cos(\omega_1 \sqrt{r-1})}\\
&=\int_1^r\frac{\cos(\omega_1\sqrt{r-z})}{\sqrt{r-z}}\hat{f_2}(\omega_1,z)\mathrm{d}z,
\end{split}
\end{equation}
a Volterra integral equation of the first kind with weakly singular kernel. 

Applying the Abel transform to both sides removes the singularity
\begin{equation}
\begin{split}
\int_1^s\frac{\hat{g}(\omega_1,r)}{\sqrt{s-r}}\mathrm{d}r&= \int_1^s\int_1^r\frac{\cos(\omega_1\sqrt{r-z})}{\sqrt{r-z}\sqrt{s-r}}\hat{f_2}(\omega_1,z)\mathrm{d}z\mathrm{d}r\\
&=\int_1^s\left[\int_z^s\frac{\cos(\omega_1\sqrt{r-z})}{\sqrt{r-z}\sqrt{s-r}}\mathrm{d}r\right]\hat{f_2}(\omega_1,z)\mathrm{d}z\\
&=\int_1^sK_1(s,z)\hat{f_2}(\omega_1,z)\mathrm{d}z,
\end{split}
\end{equation}
where
\begin{equation}
K_1(s,z)=\int_z^s\frac{\cos(\omega_1\sqrt{r-z})}{\sqrt{r-z}\sqrt{s-r}}\mathrm{d}r=\int_0^1\frac{\cos(\omega_1\sqrt{u}\sqrt{s-z})}{\sqrt{u}\sqrt{1-u}}\mathrm{d}u,
\end{equation}
after making the substitution $r=z+(s-z)u$. Hence $K_1(s,s)=\pi$ and the first derivative with respect to $s$ is
\begin{equation}
\begin{split}
\frac{\mathrm{d}}{\mathrm{d}s}K_1(s,z)&=-\omega_1^2\int_0^1\frac{u}{\sqrt{u}\sqrt{1-u}}\cdot \frac{\sin(\omega_1\sqrt{u}\sqrt{s-z})}{\omega_1\sqrt{u}\sqrt{s-z}}\mathrm{d}u\\
&=-\omega_1^2\int_0^1\frac{\sqrt{u}}{\sqrt{1-u}}\cdot \sinc(\omega_1\sqrt{u}\sqrt{s-z})\mathrm{d}u.
\end{split}
\end{equation}
It follows that
\begin{equation}
\begin{split}
\hat{g_1}(\omega_1,s)&=\frac{1}{\pi}\frac{\mathrm{d}}{\mathrm{d}s}\int_1^s\frac{\hat{g}(\omega_1,r)}{\sqrt{s-r}}\mathrm{d}r\\
&=-\frac{\omega_1^2}{\pi}\int_1^s\int_0^1\frac{\sqrt{u}}{\sqrt{1-u}}\cdot \sinc(\omega_1\sqrt{u}\sqrt{s-z})\mathrm{d}u\hat{f_2}(\omega_1,z)\mathrm{d}z+ \hat{f_2}(\omega_1,s)\\
&=-\frac{\omega_1^2}{\pi}\int_1^sK_2(s,z)\hat{f_2}(\omega_1,z)\mathrm{d}z+\hat{f_2}(\omega_1,s),
\end{split}
\end{equation}
with $|K_2|<\frac{\pi}{2}$, a bounded kernel. The Volterra equation of the second kind above can now be solved explicitly by Theorem \ref{volt}, and hence $\hat{f_2}(\omega_1,z)$ is uniquely determined by $\mathcal{T}f$ for $|\omega_1|<\frac{\pi}{2\sqrt{r_m-1}}$ and $z\in [1,r^2_m]$ (and hence for all $z\in\mathbb{R}$ due to the support restrictions on $f$). Hence we have a unique determination of the 2-D Fourier transform of $f_2$ on the open band $B=\{|\omega_1|<\frac{\pi}{2\sqrt{r_m-1}},\omega_3\in\mathbb{R}\}$, where $\omega_3$ is dual to $z$. As $f_2$ is compactly supported, its Fourier transform is analytic by the Paley-Weiner-Schwartz theorem, and hence $\hat{f_2}$ is uniquely determined everywhere on the plane by analytic continuation from $B$. It follows that $\mathcal{T}$ is injective. 
\end{proof}
\end{theorem}
\begin{remark}
The inversion process presented above uses analytic continuation in the Fourier domain to recover the density uniquely. The inversion is not stable however (due to analytic continuation) and severely ill-posed, in the sense that the solution is not bounded in any Sobolev space. In \cite{2D8} a toric section transform is considered in a rotational geometry. A toric section transform may be written as the sum of two circular Radon transforms (as in \cite{2D8}, or equivalently here as the sum of four semicircle transforms). When the circle transforms are considered separately, the injectivity follows from the results of \cite{quinto}, and the inversion is stable (as in \cite{2D8}). When the sum is considered however there exist image artefacts in the reconstruction (proven microlocally in \cite{2D8} in rotational geometries). In our case we require additional intuition to prove injectivity, since we are considering the sum of two circle transforms with circle centers on a line. We wonder if similar artefacts to those of \cite{2D8} may be present also in translational geometries for the toric section transform.
\end{remark}
\section{The three dimensional case}
\label{3D}
Let the apple $A(r)=\cup_{j=1}^2A_j(r)$ (as in figure \ref{fig1}) be written as the disjoint union of the surfaces of revolution of two semicircles, which are parameterized by
\begin{equation}
x_1=(R+\sqrt{r^2-(z-2)^2})\cos\varphi,\ \ \ y_1=(R+\sqrt{r^2-(z-2)^2})\sin\varphi\ \ \ (\text{for $A_1(r)$})
\end{equation}
and
\begin{equation}
x_2=(R-\sqrt{r^2-(z-2)^2})\cos\varphi,\ \ \ y_2=(R-\sqrt{r^2-(z-2)^2})\sin\varphi \ \ \ (\text{for $A_2(r)$})
\end{equation}
for $2-r<z<1$ and $\varphi\in[0,2\pi]$. Together, the above parameterizations describe the set of points on an apple surface (a spindle torus with the lemon part removed). Let $f\in L^2_0(\Omega)$ for $\Omega$ compactly supported in $\{2-r_m<z<1\}\subset\mathbb{R}^3$. We define the apple transform $\mathcal{A}: L^2_0(\Omega)\to \mathcal{A}(L^2_0(\Omega))$ in the translational geometry as
\begin{equation}
\begin{split}
\mathcal{A}f(x_0,y_0,r)=\int_{T_{x_0,y_0}(A(r))}f\mathrm{d}A=\sum_{j=1}^2\int_{T_{x_0,y_0}(A_j(r))}f\mathrm{d}A_j,
\end{split}
\end{equation}
where $\mathrm{d}A_j$ denotes the surface area measure on $A_j(r)$ and $T_{x_0,y_0}(x,y,z)=(x+x_0,y+y_0,z)$ denotes a translation in the $xy$ plane to the point $(x_0,y_0)$. We now proceed in a similar vein to the 2-D case.
\begin{proposition}
Let $f\in L^2_0(\Omega)$ for $\Omega$ compactly supported in $\{2-r_m<z<1\}\subset\mathbb{R}^3$, and let $f_1(x,y,z)=f(x,y,2-z)$. Then
\begin{equation}
\label{Aform}
\begin{split}
\mathcal{A}f(x_0,y_0,r)=\int_{-\pi}^{\pi}\int_1^r\frac{r}{\sqrt{r^2-z^2}}\sum_{j=1}^2\rho_j f_1(\rho_j\cos\varphi+x_0,\rho_j\sin\varphi+y_0,z)\mid_{\rho_j=R+(-1)^{j}\sqrt{r^2-z^2}}\mathrm{d}z\mathrm{d}\varphi,
\end{split}
\end{equation}
\begin{proof}
Let $\mathrm{d}s_j$, for $j=1,2 $, be the circular arc measures as in equation \eqref{arc}. Then the surface measures $\mathrm{d}A_j$ for $j=1,2$ are given by \cite[page 4]{3D2}
\begin{equation}
\mathrm{d}A_j=\mathrm{d}s_j\rho_j\mathrm{d}z\mathrm{d}\varphi=\frac{r}{\sqrt{r^2-(z-2)^2}}(R+(-1)^{j}\sqrt{r^2-(z-2)^2})\mathrm{d}z\mathrm{d}\varphi.
\end{equation}
It follows that
\begin{equation}
\begin{split}
\mathcal{A}f(x_0,y_0,r)&=\sum_{j=1}^2\int_{T_{x_0,y_0}(A_j(r))}f\mathrm{d}A_j\\
&=\sum_{j=1}^2\int_{-\pi}^{\pi}\int_h^1\frac{r}{\sqrt{r^2-(z-2)^2}}(R+(-1)^{j}\sqrt{r^2-(z-2)^2}) f(x_j+x_0,y_j+y_0,z)\mathrm{d}z\mathrm{d}\varphi\\
&=\int_{-\pi}^{\pi}\int_1^r\frac{r}{\sqrt{r^2-z^2}}\sum_{j=1}^2\rho_j f_1(\rho_j\cos\varphi+x_0,\rho_j\sin\varphi+y_0,z)\mid_{\rho_j=R+(-1)^{j}\sqrt{r^2-z^2}}\mathrm{d}z\mathrm{d}\varphi
\end{split}
\end{equation}
\end{proof}
\end{proposition}
\begin{proposition}
\label{prp3}
Let $r_m>1$ and let $\Omega\subset\{2-r_m<z<1\}\subset\mathbb{R}^3$ be compact. Then the image $\mathcal{A}(L^2_0(\Omega))\subset L^2([1,r_m]\times\mathbb{R}^2)$.
\begin{proof}
Taking the Fourier transform in the $x_0$ and $y_0$ variables of \eqref{Aform} yields
\begin{equation}
\begin{split}
\widehat{\mathcal{A}f}(\omega_1,\omega_2,r)&=\int_{\mathbb{R}^2}\mathcal{A}f(x_0,y_0,r)e^{-i(x_0,y_0)\cdot(\omega_1,\omega_2)}\mathrm{d}x_0\mathrm{d}y_0\\
&=\int_1^r\frac{rK(r,z)}{\sqrt{r^2-z^2}}\hat{f_1}(\omega_1,\omega_2,z)\mathrm{d}z,
\end{split}
\end{equation}
where $\omega_2$ is dual to $y$ (or $y_0$) and
\begin{equation}
\hspace{-2cm}
\begin{split}
K(r,z)&=\sum_{j=1}^2\rho_j \int_{-\pi}^{\pi}\exp(-i\rho_j(\omega_1\cos\varphi+\omega_2\sin\varphi))\mid_{\rho_j=R+(-1)^{j}\sqrt{r^2-z^2}}\mathrm{d}\varphi
\end{split}
\end{equation}
and $|K|<M$ is bounded. We have
\begin{equation}
\begin{split}
\|\mathcal{A}f\|^2_{L^2([1,r_m]\times\mathbb{R}^2)}&=\|\widehat{\mathcal{A}f}\|^2_{L^2([1,r_m]\times\mathbb{R}^2)}\ \ \ \text{by the Plancharel theorem}\\
&=\int_{-\infty}^{\infty}\int_{-\infty}^{\infty}\int_1^{r_m}\left(\int_1^r\frac{rK(r,z)}{\sqrt{r^2-z^2}}\hat{f_1}(\omega_1,\omega_2,z)\mathrm{d}z\right)^2\mathrm{d}r\mathrm{d}\omega_1\mathrm{d}\omega_2\\
&\leq M^2r^2_m\int_{-\infty}^{\infty}\int_{-\infty}^{\infty}\int_1^{r_m}\left(\int_1^r\frac{\hat{f_1}(\omega_1,\omega_2,z)}{\sqrt{r^2-z^2}}\mathrm{d}z\right)^2\mathrm{d}r\mathrm{d}\omega_1\mathrm{d}\omega_2\\
&\leq \frac{(Mr_m)^2}{2}\int_{-\infty}^{\infty}\int_{-\infty}^{\infty}\int_1^{r_m}\left(\int_1^r\hat{f_1}(\omega_1,\omega_2,z)L(r,z)\mathrm{d}z\right)^2\mathrm{d}r\mathrm{d}\omega_1\mathrm{d}\omega_2\\
\end{split}
\end{equation}
where $L(r,z)=(r-z)^{1/2}$. From here the proof follows the same arguments as in the proof of Proposition \ref{prp2}.
\end{proof}
\end{proposition}
We now have our second main theorem, which proves the injectivity of the apple transform on $L^2_0(\Omega)$.
\begin{theorem}
\label{main2}
Let $r_m>1$, let $0<\delta<r_m-1$ and let $\Omega\subset\{2-r_m<z<1-\delta\}\subset\mathbb{R}^3$ be compact and bounded away from $\{z=1\}$. Then the apple transform $\mathcal{A}: L^2_0(\Omega)\to L^2([1,r_m]\times\mathbb{R}^2)$ is injective.
\begin{proof}
Let
$$(\omega_1,\omega_2)=\sqrt{\omega_1^2+\omega_2^2}(\cos\varphi_{\omega},\sin\varphi_{\omega})=|\omega|(\cos\varphi_{\omega},\sin\varphi_{\omega})$$
Then from Proposition \ref{prp3} we have
\begin{equation}
\begin{split}
\widehat{\mathcal{A}f}(\omega_1,\omega_2,r)=\int_1^r\frac{rK(r,z)}{\sqrt{r^2-z^2}}\hat{f_1}(\omega_1,\omega_2,z)\mathrm{d}z,
\end{split}
\end{equation}
where
\begin{equation}
\hspace{-2cm}
\begin{split}
K(r,z)&=\sum_{j=1}^2\rho_j \int_{-\pi}^{\pi}\exp(-i\rho_j(\omega_1\cos\varphi+\omega_2\sin\varphi))\mid_{\rho_j=R+(-1)^{j}\sqrt{r^2-z^2}}\mathrm{d}\varphi\\
&=\sum_{j=1}^2\rho_j \int_{-\pi}^{\pi}\exp(i\rho_j|\omega|\cos(\varphi-\varphi_{\omega}-\pi))\mid_{\rho_j=R+(-1)^{j}\sqrt{r^2-z^2}}\mathrm{d}\varphi\\
&=\sum_{j=1}^2\rho_j \int_{-\pi}^{\pi}\exp(i\rho_j|\omega|\cos\varphi)\mid_{\rho_j=R+(-1)^{j}\sqrt{r^2-z^2}}\mathrm{d}\varphi,\ \ \ \text{(due to periodicity)}\\
&=2\pi\sum_{j=1}^2(R+(-1)^{j}\sqrt{r^2-z^2}) J_0(|\omega|(R+(-1)^{j}\sqrt{r^2-z^2})),
\end{split}
\end{equation}
where $J_0$ is a Bessel function of the first kind of order zero. Letting $\hat{f_2}(\omega_1,\omega_2,z)=\frac{\hat{f_1}(\omega_1,\omega_2,\sqrt{z})}{2\sqrt{z}}$, we have
\begin{equation}
\begin{split}
\hat{g}(\omega_1,\omega_2,r)&=\frac{\widehat{\mathcal{A}f}(\omega_1,\omega_2,\sqrt{r})}{\sqrt{r}}\\
&=\int_1^r\frac{K(\sqrt{r},\sqrt{z})}{\sqrt{r-z}}\hat{f_2}(\omega_1,\omega_2,z)\mathrm{d}z.
\end{split}
\end{equation}
It follows that
\begin{equation}
\label{17}
\begin{split}
\int_1^s\frac{\hat{g}(\omega_1,\omega_2,r)}{\sqrt{s-r}}\mathrm{d}r&= \int_1^s\int_1^r\frac{K(\sqrt{r},\sqrt{z})}{\sqrt{r-z}\sqrt{s-r}}\hat{f_2}(\omega_1,\omega_2,z)\mathrm{d}z\mathrm{d}r\\
&=\int_1^s\left[\int_z^s\frac{K(\sqrt{r},\sqrt{z})}{\sqrt{r-z}\sqrt{s-r}}\mathrm{d}r\right]\hat{f_2}(\omega_1,\omega_2,z)\mathrm{d}z\\
&=\int_1^sK_1(s,z)\hat{f_2}(\omega_1,\omega_2,z)\mathrm{d}z,
\end{split}
\end{equation}
where
\begin{equation}
\begin{split}
K_1(s,z)&=\int_z^s\frac{K(\sqrt{r},\sqrt{z})}{\sqrt{r-z}\sqrt{s-r}}\mathrm{d}r\\
&=2\pi\int_0^1\frac{\sum_{j=1}^2(R+(-1)^{j}\sqrt{s-z}\sqrt{u}) J_0(|\omega|(R+(-1)^{j}\sqrt{s-z}\sqrt{u}))}{\sqrt{u}\sqrt{1-u}}\mathrm{d}u\\
&=\int_0^1\frac{H(s,z,u)}{\sqrt{u}\sqrt{1-u}}\mathrm{d}u,
\end{split}
\end{equation}
where $R=\sqrt{(z-1)+(s-z)u}$. To calculate the derivatives of $H$ we use the previous expression for the Bessel function
\begin{equation}
J_0(\nu)=\frac{1}{2\pi}\int_{-\pi}^{\pi}e^{i\nu\cos\varphi}\mathrm{d}\varphi.
\end{equation}
Hence
\begin{equation}
\hspace{-2cm}
\begin{split}
H(s,z,u)&=2R\int_{-\pi}^{\pi}e^{i|\omega|R\cos\varphi}\cos(|\omega|\sqrt{u}\sqrt{s-z}\cos\varphi)\mathrm{d}\varphi\\
&+2i\sqrt{u}\sqrt{s-z}\int_{-\pi}^{\pi}e^{i|\omega|R\cos\varphi}\sin(|\omega|\sqrt{u}\sqrt{s-z}\cos\varphi)\mathrm{d}\varphi\\
&=\int_{-\pi}^{\pi}h(s,z,u,\varphi)(h_1(s,z,u,\varphi)+ih_2(s,z,u,\varphi))\mathrm{d}\varphi,
\end{split}
\end{equation}
where $h(s,z,u,\varphi)=e^{i|\omega|R\cos\varphi}$,
\begin{equation}
h_1(s,z,u,\varphi)=2R\cos(|\omega|\sqrt{u}\sqrt{s-z}\cos\varphi)
\end{equation}
and
\begin{equation}
h_2(s,z,u,\varphi)=2\sqrt{u}\sqrt{s-z}\sin(|\omega|\sqrt{u}\sqrt{s-z}\cos\varphi).
\end{equation}
The first partial derivative of $h$ with respect to $s$ is
\begin{equation}
\frac{\mathrm{d}}{\mathrm{d}s}h(s,z,u,\varphi)=\frac{i|\omega|u\cos\varphi}{2R}e^{i|\omega|R\cos\varphi},
\end{equation}
which is bounded on the support of $f$, since $\text{supp}(f)$ is bounded away from $\{z=1\}$ by assumption (and hence $R$ is bounded away from 0). We have
\begin{equation}
\begin{split}
\frac{\mathrm{d}}{\mathrm{d}s}h_1(s,z,u,\varphi)&=\frac{u}{R}\cos(|\omega|\sqrt{u}\sqrt{s-z}\cos\varphi)-R\frac{|\omega|\sqrt{u}\cos\varphi\sin(|\omega|\sqrt{u}\sqrt{s-z}\cos\varphi)}{\sqrt{s-z}}\\
&=\frac{u}{R}\cos(|\omega|\sqrt{u}\sqrt{s-z}\cos\varphi)-R|\omega|^2u\cos^2\varphi\sinc(|\omega|\sqrt{u}\sqrt{s-z}\cos\varphi),
\end{split}
\end{equation}
and
\begin{equation}
\begin{split}
\frac{\mathrm{d}}{\mathrm{d}s}h_2(s,z,u,\varphi)&=\frac{\sqrt{u}}{\sqrt{s-z}}\sin(|\omega|\sqrt{u}\sqrt{s-z}\cos\varphi)+\frac{u\sqrt{s-z}|\omega|\cos\varphi\cos(|\omega|\sqrt{u}\sqrt{s-z}\cos\varphi)}{\sqrt{s-z}}\\
&=|\omega|u\cos\varphi\sinc(|\omega|\sqrt{u}\sqrt{s-z}\cos\varphi)+u|\omega|\cos\varphi\cos(|\omega|\sqrt{u}\sqrt{s-z}\cos\varphi).
\end{split}
\end{equation}
It can now be seen from the above and by an application of the product rule, that $H_1(s,z,u)=\frac{\mathrm{d}}{\mathrm{d}s}H(s,z,u)$ is bounded on the limits of integration, where $f$ is supported. That is for $(1+\delta)^2<s<r^2_m$, $(1+\delta)^2<z<s$ and $0<u<1$, for any fixed $|\omega|$. 

Now (\ref{17}) becomes
\begin{equation}
\label{18}
\begin{split}
\hat{g_1}(\omega_1,\omega_2,r)&=\frac{\mathrm{d}}{\mathrm{d}s}\int_1^s\frac{\hat{g}(\omega_1,\omega_2,r)}{\sqrt{s-r}}\mathrm{d}r\\
&=\int_1^s\frac{\mathrm{d}}{\mathrm{d}s}K_1(s,z)\hat{f_2}(\omega_1,\omega_2,z)\mathrm{d}z+2\sqrt{s-1}J_0(|\omega|\sqrt{s-1})\hat{f_2}(\omega_1,\omega_2,s).
\end{split}
\end{equation}
Let $t_0$ denote the first root of $J_0$. Then for $|\omega|<\frac{t_0}{\sqrt{r^2_m-1}}$ and $s<r^2_m$, $J_0(|\omega|\sqrt{s-1})>0$ and we have
\begin{equation}
\label{18}
\begin{split}
\frac{\hat{g_1}(\omega_1,\omega_2,s)}{2\sqrt{s-1}J_0(|\omega|\sqrt{s-1})}&=\frac{\mathrm{d}}{\mathrm{d}s}\int_1^s\frac{\hat{g}(\omega_1,\omega_2,r)}{\sqrt{s-r}}\mathrm{d}r\\
&=\int_1^s\frac{\frac{\mathrm{d}}{\mathrm{d}s}K_1(s,z)}{2\sqrt{s-1}J_0(|\omega|\sqrt{s-1})}\hat{f_2}(\omega_1,\omega_2,z)\mathrm{d}z+\hat{f_2}(\omega_1,\omega_2,s),
\end{split}
\end{equation}
a Volterra equation of the second kind with $|\frac{\mathrm{d}}{\mathrm{d}s}K_1(s,z)|\leq\int_0^1\frac{|\frac{\mathrm{d}}{\mathrm{d}s}H(s,z,u)|}{\sqrt{u}\sqrt{1-u}}\mathrm{d}u<M$ bounded kernel on $\{(1+\delta)^2<s<r^2_m, (1+\delta)^2<z<s\}$. As in the 2-D case, and by Theorem \ref{volt}, we can reconstruct the 3-D Fourier transform of $\hat{f_2}$ uniquely on the open set 
$$B=\left\{\sqrt{\omega_1^2+\omega_2^2}<\frac{t_0}{\sqrt{r^2_m-1}}, \omega_3\in\mathbb{R}\right\}$$
and hence for all $(\omega_1,\omega_2,\omega_3)\in\mathbb{R}^3$ by analytic continuation and the Paley-Weiner-Schwartz theorem. We conclude that $\mathcal{A}$ is injective.
\end{proof}
\end{theorem}
\begin{discussion}
Theorems \ref{main1} and \ref{main2} explain how we can recover $\hat{f}$ on an open set in the Fourier domain, by bounding the limits on cosine and Bessel kernels so as to obtain an explicit solution via repeated application of Volterra operators. The set $S$ of $\omega$ for which $\hat{f}$ is known (without the need for analytic continuation) is bounded by the curves
$$c_1=\left\{\omega_1=-\frac{t_0}{\sqrt{\omega_3^2-1}}\right\},\ \ \ c_2=\left\{\omega_1=\frac{t_0}{\sqrt{\omega_3^2-1}}\right\},\ \ \ c_3=\{\omega_3=1\}$$
in 2-D and by the surface of revolution of $c_1$ (or $c_2$ due to symmetry) about $\omega_3$ and the plane $\{\omega_3=1\}$ in 3-D (this is evident from the proofs of theorems \ref{main1} and \ref{main2}). See figure \ref{fig2}. While the entirety of the Fourier space is uniquely determined from toric section (2-D) and apple (3-D) integral data, the reconstruction outside of $S$ is not stable using the proposed inversion formulae, due to analytic continuation. Further, we have greater recovery of $\hat{f}$ in the Fourier space (without the need for analytic continuation) as $\omega_3\to1$ (and at a rate $1/\omega_3$), and hence we can expect an increasingly stable recovery of $f$ as $z\to 1$ (for points in the scanning region closer to the detector line and the photon source). 
\begin{figure}[!h]
\centering
\begin{tikzpicture}[scale=4]
\draw [very thick] (-1.5,1-0.5)--(1.5,1-0.5)node[right] {detector line};
\draw [very thick] (-1.5,1.5-0.5)--(1.5,1.5-0.5) node[right] {source line};
\draw [->] (0,0.5-0.5)--(0,1.7-0.5)node[right] {$\omega_3$};
\draw [->] (0,0.5-0.5)--(1,0.5-0.5)node[below] {$\omega_1$};
\node at (-0.05,-0.05) {$O$};
\draw[scale=0.5,domain=0.3:3,smooth,variable=\x,blue] plot ({1/\x},{-\x+1});
\draw[scale=0.5,domain=0.3:3,smooth,variable=\x,blue] plot ({-1/\x},{-\x+1});
\draw[blue] (-0.1,-1)--(0.1,-1);
\draw[blue] (-0.15,-0.8)--(0.15,-0.8);
\draw[blue] (-0.2,-0.6)--(0.2,-0.6);
\draw[blue] (-0.25,-0.4)--(0.25,-0.4);
\draw[blue] (-0.3,-0.2)--(0.3,-0.2);
\draw[blue] (-0.45,0.05)--(0.45,0.05);
\draw[blue] (-0.85,0.25)--(0.85,0.25);
\draw[blue] (-1.5,0.45)--(1.5,0.45);
\node at (-0.68,-0.2) {$\omega_1=-\frac{t_0}{\sqrt{\omega_3^2-1}}$};
\node at (0.68,-0.2) {$\omega_1=\frac{t_0}{\sqrt{\omega_3^2-1}}$};
\end{tikzpicture}
\caption{The RTT geometry (with the same dimensions as in figure \ref{fig1} and a visualisation of the set
$S=\{\omega_1>-\frac{t_0}{\sqrt{\omega_3^2-1}}\}\cap\{\omega_1<\frac{t_0}{\sqrt{\omega_3^2-1}}\}\cap\{\omega_3<1\}$ (shaded in blue) of $\omega$ for which $\hat{f}(\omega)$ is determined without the need for analytic continuation. In 2-D $t_0=\frac{\pi}{2}$. In 3-D $t_0$ is the first root of $J_0$ and the set of $\omega$ recovered without analytic continuation is the volume of revolution of $S$ about $z$.}
\label{fig2}
\end{figure}
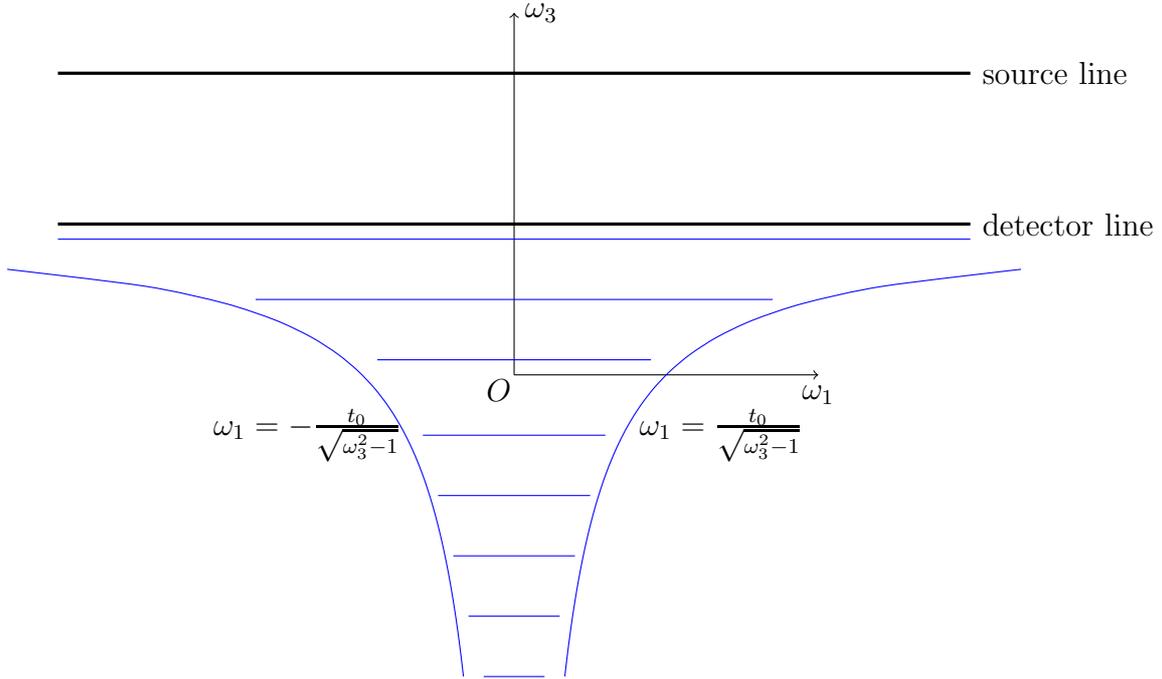
\end{discussion}

\section{A generalization to the surfaces of revolution of $C^1$ curves}
\label{gen}
So far we have considered the surfaces of revolution of semicircles in translational geometries and proven injectivity and explicit invertibility results for classes of $L^2$ functions of compact support. Here we consider a more general class of surfaces, which are the surfaces of revolution of $C^1$ curves, and we prove injectivity results on the set of continuous functions of compact support. 

For some $r_m>1$, let $\rho_j(r,z)\in C^1([1,r_m]^2)$, for $j=1,\ldots,m$, parameterize a finite set of surfaces of revolution in cartesian coordinates
\begin{equation}
x_j=\rho_j(r,z)\cos\varphi,\ \ \ y_j=\rho_j(r,z)\sin\varphi, \ \ \ 1<z<r_m,\ \ \  0\leq\varphi\leq 2\pi.
\end{equation}
Let $\Omega\subset\{2-r_m<z<1\}\subset\mathbb{R}^3$ be compact. Then we define the generalized apple transform $R_A : C_0(\Omega)\to C([1,r_m]\times\mathbb{R}^2)$
\begin{equation}
\label{Grad}
R_Af(x_0,y_0,r)=\sum_{j=1}^m\int_{-\pi}^{\pi}\int_1^r\sqrt{1+\left(\frac{\mathrm{d}\rho_j}{\mathrm{d}z}\right)^2}\rho_j(r,z)f(\rho_j(r,z)\cos\varphi+x_0,\rho_j(r,z)\sin\varphi+y_0,2-z)\mathrm{d}z\mathrm{d}\varphi
\end{equation}
We now have our third main theorem, which is a natural extension of the results of sections \ref{2D} and \ref{3D}.
\begin{theorem}
Let $\Omega\subset\{2-r_m<z<1\}\subset\mathbb{R}^3$ be compact and let $\rho_j(r,z)\in C^1([1,r_m]^2)$, for $j=1,\ldots,m$ be a finite set of $C^1$ curves. Then the generalized apple transform $R_A : C_0(\Omega)\to C([1,r_m]\times\mathbb{R}^2)$ in \eqref{Grad} is injective.
\begin{proof}
We have
\begin{equation}
\begin{split}
\widehat{R_Af}(\omega_1,\omega_2,r)=\int_1^r\sqrt{1+\left(\frac{\mathrm{d}\rho_j}{\mathrm{d}z}\right)^2}K(r,z)\hat{f_1}(\omega_1,\omega_2,z)\mathrm{d}z,
\end{split}
\end{equation}
where $f_1(x,y,z)=f(x,y,2-z)$ as before, and
\begin{equation}
\hspace{-2cm}
\begin{split}
K(r,z)&=\sum_{j=1}^m\rho_j(r,z) \int_{-\pi}^{\pi}\exp(-i\rho_j(r,z)(\omega_1\cos\varphi+\omega_2\sin\varphi)\mathrm{d}\varphi\\
&=2\pi\sum_{j=1}^m \rho_j(r,z)J_0(|\omega|\rho_j(r,z)),
\end{split}
\end{equation}
where $J_0$ is a Bessel function of the first kind and $|\omega|=\sqrt{\omega_1^2+\omega_2^2}$. Let $\rho_j<M_j$ and without loss of generality we can assume that $\rho_j\geq 0$ for all $j$. Now let
\begin{equation}
\label{corm}
\begin{split}
0=\int_1^r\sqrt{1+\left(\frac{\mathrm{d}\rho_j}{\mathrm{d}z}\right)^2}K(r,z)\hat{f_1}(\omega_1,\omega_2,z)\mathrm{d}z,
\end{split}
\end{equation}
for $r\in[1,r_m]$. Let $t_0$ be the first root of $J_0$ and let $|\omega|<\frac{t_0}{M}$, where $M=\max_j M_j$. Then for such $|w|$, the integrand \eqref{corm} is positive except for $\hat{f_1}$, and hence the only continuous solution to \eqref{corm} is $\hat{f_1}(\omega_1,\omega_2,z)=0$ for $z\in[1,r_m]$ and $|\omega|<\frac{t_0}{M}$. Hence it follows that the 3-D Fourier transform of $f_1$ is zero on the open set
$$B=\left\{\sqrt{\omega_1^2+\omega_2^2}<\frac{t_0}{M}, \omega_3\in\mathbb{R}\right\}.$$
The result follows by the Paley-Weiner-Schwartz theorem and analytic continuation.
\end{proof}
\end{theorem}
\begin{remark}
The above proof uses the same arguments to that of Cormack \cite[page 2724]{cor} in proofs of injectivity for Volterra integral equations, which suggest that Volterra operators of the first kind with positive kernel, are injective on the domain of continuous functions. We wonder if such ideas could be extended to prove injectivity in the $L^2$ case.
\end{remark}
\section{Conclusions and further work}
Here we have presented new injectivity results and explicit inversion formulae for CST problems in translational geometries. We considered the problem of electron density reconstruction from sets of toric section (in 2-D) and apple integral data (in 3-D). In section \ref{2D} we introduced a new two dimensional toric section Radon transform $\mathcal{T}$ which describes the integrals of an $L^2$ density over the set of toric sections translated along a line, whose central axis is vertical (parallel to the $z$ axis). Here we provided boundedness theorems for $\mathcal{T}$ in $L^2$ and went on the prove the injectivity and explicit invertibility of $\mathcal{T}$ on $L^2_0$. After proving the $L^2$ injectivity in the 2-D case, we considered the 3-D case in section \ref{3D}, where we introduced a new apple Radon transform $\mathcal{A}$ in a three dimensional scanning modality previously introduced in \cite{3D2} (configuration (d) on page 5). The transformation $\mathcal{A}$ takes an $L^2$ density to its sinogram of integrals over vertical apples which are translated in the $xy$ plane. In a similar vein to the 2-D case, we showed the continuity of the apple operator in $L^2$ and then went on to prove the injectivity and explicit invertibility of $\mathcal{A}$ in $L^2_0$. The proofs of injectivity in the 2-D and 3-D cases followed similar ideas, solving a set of 1-D Volterra operators to recover the density in an open subset of the Fourier domain, and then extending to the entire Fourier space uniquely by analytic continuation (which was possible using Paley-Weiner-Schwartz ideas as the densities considered were of compact support). In the final part of the paper in section \ref{gen} we presented a generalized apple Radon transform $R_A$, which describes the integrals of a density over the surfaces of revolution of $C^1$ curves in the translational geometry. Here we proved the injectivity of $R_A$ on the space of continuous functions of compact support using the results of Cormack \cite{cor}.

The work of \cite{2D8} characterizes microlocally the image artefacts in a reconstruction from toric section integral data in rotational geometries. Further, an iterative reconstruction scheme using TV is found to be effective in suppressing the discovered artefacts. In further work we aim to consider the microlocal properties of $\mathcal{T}$ and $\mathcal{A}$ to investigate the existence of image artefacts. The inversion formulae presented here rely on analytic continuation in the Fourier domain, and hence the inversion process is severely ill-posed, and we would expect a large amplification of the measurement error in the solution. However, this is only one possible inversion approach, and is not proof that the problem of recovering $f$ from $\mathcal{T}f$ or $\mathcal{A}f$ is severely ill-posed. We aim to clarify this in future work. Following this, we aim to devise a reconstruction method as in \cite{2D8} which can offer a practically useful image quality for the desired application in airport baggage screening.
\section*{Acknowledgements}
We would like to thank Prof. Eric Todd Quinto for his helpful discussions, comments and insight towards the results presented in this paper. 

This material is based upon work supported by the U.S. Department of Homeland Security, Science and Technology Directorate, Office of University Programs, under Grant Award 2013-ST-061-ED0001. The views and conclusions contained in this document are those of the authors and should not be interpreted as necessarily  epresenting the official policies, either expressed or implied, of the U.S. Department of Homeland Security.

\end{document}